\documentclass[10 pt]{amsart}

\usepackage{amsmath}
\usepackage{amsfonts}
\usepackage{amssymb}
\usepackage{amstext}
\usepackage{amsbsy}
\usepackage{amsopn}
\usepackage{amsthm}
\usepackage{amsxtra}
\usepackage{graphicx}
\usepackage{hyperref}
\usepackage{color}
\usepackage[all]{xy}
\usepackage{enumitem}
\usepackage{tikz-cd}
\usepackage{mathtools}
\usepackage{array}
\usepackage{comment}
\usepackage{bm}
\usepackage{soul}

\newtheorem{theorem}{Theorem}[section]
\newtheorem{lemma}[theorem]{Lemma}
\newtheorem{proposition}[theorem]{Proposition}
\newtheorem{corollary}[theorem]{Corollary}

\theoremstyle{definition}
\newtheorem*{definition*}{Definition}

\theoremstyle{remark}

\newtheorem*{remark*}{Remark}

\newcommand{\CA}{\mathcal{A}}

\newcommand{\CM}{\mathcal{M}}

\newcommand{\CK}{\mathcal{K}}

\newcommand{\CF}{\mathcal{F}}

\newcommand{\CW}{\mathcal{W}}
\newcommand{\CE}{\mathcal{E}}
\newcommand{\CY}{\mathcal{Y}}
\newcommand{\CX}{\mathcal{X}}
\newcommand{\CD}{\mathcal{D}}

\newcommand{\inv}[1]{\textnormal{Inv}(\sigma_{n})}
\newcommand{\infinv}[1]{\textnormal{Inv}^{\infty}(\sigma_{n})}

\newcommand{\T}{{\mathbb T}}

\newcommand{\dimH}{\dim_{\mathrm H}}
\newcommand{\htop}{h_{\mathrm{top}}}

\usepackage{todonotes}
\begin{document}

\title{Invariant random compacts for $\times p, \times q$ on the circle}

\author{Scott Schmieding}
\address{Department of Mathematics, 
107 McAllister Building
University Park, State College, PA 16802}
\email{sks7247@psu.edu}

\thanks{SS was partially supported by NSF grant DMS-2247553.}

\begin{abstract}
We prove that for multiplicatively independent natural numbers $p,q \ge 2$, if $\mu$ is a Borel probability measure on the space $\CK(\T)$ of nonempty compact subsets of the circle $\T$ which is invariant under the induced actions of both multiplication by $p$ and $q$, and $\mu$-a.e. $Y$ has positive Hausdorff dimension, then $\mu$-a.e. $Y$ is equal to $\T$.
\end{abstract}
\maketitle

\section{Introduction}
Let $\T=\mathbb{R}/\mathbb{Z}$ denote the circle and for $n \ge 2$ let $T_{n} \colon \T \to \T$ by multiplication by $n$ defined by $T_{n}(x) = nx \textrm{ mod } 1$. Natural numbers $p,q$ are said to be multiplicatively independent if they are not both powers of a common integer. In~\cite{Furstenberg1967}, Furstenberg proved that if $Y \subset \T$ is compact, infinite, and invariant under both $T_{p}$ and $T_{q}$ for $p,q \ge 2$ multiplicatively independent, then $Y = \T$. 


In~\cite{KSirc}, Kra and the author introduced the notion of an invariant random compact for an action. Given a semigroup $G$ acting on a compact metric space $X$, there is an induced action on the space $\CK(X)$ of nonempty compact subsets of $X$, considered with the Hausdorff metric. An \emph{invariant random compact} (IRC for short) for the action is a Borel probability measure $\mu$ on $\CK(X)$ which is invariant under the induced action of $G$ on $\CK(X)$. If $Y \subset X$ is a fully-invariant subset (meaning $gY=Y$ for every $g \in G$), then $\delta_{Y}$ is an IRC; in this sense, one may consider IRCs as a generalization of fully-invariant sets.

Of particular interest is when the IRC $\mu$ is \emph{nonfinitary}, meaning $\mu$-almost every compact set $Y$ is infinite. In~\cite{KSicrigid}, it was shown that the only nonfinitary IRC for the multiplicative semigroup $\mathbb{N}$ acting by multiplication on $\T$ is $\delta_{\T}$. 

Motivated by Furstenberg's theorem, in~\cite{KSirc}, the following question was posed.\\

\textbf{Question: } For $p,q \ge 2$ multiplicatively independent, if $\mu$ is an IRC for multiplication by $p,q$ on $\T$ and $\mu$-a.e. compact set $Y$ is infinite, must $\mu=\delta_{\T}$?\\

Here we show that, with a dimension assumption on the IRC, the answer is yes. More precisely, we prove the following theorem.

\begin{theorem}\label{thm:mainthm}
Suppose $p,q \ge 2$ are multiplicatively independent natural numbers and $\mu$ is a Borel probability measure on $\CK(\T)$ which is invariant under the induced actions of both $T_{p}$ and $T_{q}$. Suppose further that $\mu$-a.e. $Y$ has positive Hausdorff dimension. Then $\mu$-a.e. $Y$ is equal to $\T$.
\end{theorem}

We obtain from this the following corollary.
\begin{corollary}\label{cor:corollary1}
Suppose $p,q \ge 2$ are multiplicatively independent and $Y \subset \T$ is compact, nonempty, and $s$-Ahlfors-regular for some $s>0$. Then for every $\varepsilon > 0$, 
\begin{equation}
\frac{1}{N^{2}}\lvert\{0 \le a,b < N : T_{p}^{a}T_{q}^{b}(Y) \textrm{ is } \varepsilon\textrm{-dense in } \T\}\rvert \to 1 \quad \textrm{ as } N \to \infty.
\end{equation}
\end{corollary}

It follows from~\cite[Theorem 1.5]{BerendPeres} that there exists an infinite compact set $Y \subset \T$ and $\varepsilon > 0$ such that $T_{2}^{a}T_{6}^{b}(Y)$ is not $\varepsilon$-dense in $\T$ for every $a,b \ge 0$; in the language used there, this means the multiplicative semigroup generated by $2$ and $6$ is not a Glasner set for the action of multiplication on the circle. To our knowledge, whether the multiplicative semigroup generated by $2$ and $3$ is Glasner or not for the action on the circle is not known. The above corollary gives the dense dilations necessary for the Glasner property for all multiplicatively independent $p,q$, but assuming quite a strong condition on the set $Y$, namely, $s$-Ahlfors regularity for some $s>0$. 

The majority of the paper is devoted to proving Theorem~\ref{thm:mainthm}. Corollary~\ref{cor:corollary1} is proved after. 

We briefly outline the idea for the proof of Theorem~\ref{thm:mainthm}. The starting point is based on the proof of Theorem 4.2 in~\cite{KSirc}. For the proof, it is convenient to work at a symbolic level, but we sketch the idea on the circle, as the overall idea is similar.

Given an IRC $\mu$ for $T_{p},T_{q}$, consider $\mathcal{Z} = \{(Y,x) \in \mathcal{K}(\T) \times \T : x \in Y\}$. Letting $\tilde{T}_{k}$ denote the map induced by $T_{k}$ on $\CK(\T)$, for $k=p,q$ we have $S_{k} = \tilde{T}_{k} \times T_{k}$ acting on $\mathcal{Z}$, and we lift $\mu$ to a jointly $S_{p},S_{q}$-invariant measure $\eta$ on $\mathcal{Z}$. After passing to the invertible case (which we ignore notationally for this discussion), we then disintegrate $\eta = \int\eta_{Y}\, d\mu$ over the first component of $\mathcal{Z}$ and arrive at a family $\eta_{Y}$ of measures with each $\eta_{Y}$ supported on $Y$ where $Y$ is sampled via $\mu$. Averaging the self-convolution produces a jointly invariant measure $\xi$, and the idea is then to apply a theorem of Hochman-Shmerkin to conclude that $\xi$ is Lebesgue, from which the theorem shortly after follows. To apply Hochman-Shmerkin however, we need that $\xi$ has positive dimension; this is where the dimension assumption on $\mu$ comes in, and the majority of the proof is devoted to ensuring we can construct such $\xi$ with positive dimension.

\subsection{Acknowledgments}
The author is grateful to Bryna Kra for several helpful conversations related to the results here. AI was used on a draft to mathematically audit the paper, find some references, correct typos, with the proof of Lemma 3.4 slightly corrected and simplified as a result of the audit.

\section{Notation and background}
We let $\T=\mathbb{R}/\mathbb{Z}$ denote the circle. Given $n \ge 2$ we define $T_{n} \colon \T \to \T$ by $T_{n}(x) = nx \textrm{ mod } 1$. We write $\lambda_{\T}$ for Lebesgue measure on $\T$. Logarithms will always be taken base two. For a space $X$, we write $\mathcal{M}(X)$ for the space of Borel probability measures on $X$.

We consider the full shift over $\mathbb{Z}$ on $n$ symbols defined by $X_{n} = \{0,\ldots,n-1\}^{\mathbb{Z}}$ and let $\sigma_{n} \colon X_{n} \to X_{n}$ be the left shift map. We also consider the one-sided full shifts $X_{n}^{+} = \{0,\ldots,n-1\}^{\mathbb{N}_{0}}$ and let $\sigma_{n}^{+} \colon X_{n}^{+} \to X_{n}^{+}$ denote the respective left shift map. On both $X_{n}$ and $X_{n}^{+}$ we use the metric
$$d(x_{k},y_{k}) = 2^{-\inf \{|k| : x_{k} \ne y_{k}\}}.$$

There is a factor map
\[
P_{n} \colon X_{n} \to X_{n}^{+}
\]
\[
P_{n}((x_{k})) = (x_{k})_{k \in \mathbb{N}_{0}}.
\]

Let $\lambda_{n}^{+}$ denote Haar measure on $X_{n}^{+}$ and $\lambda_{n}$ denote Haar measure on $X_{n}$; equivalently, $\lambda_{n}^{+}$ and $\lambda_{n}$ are the $(\frac{1}{n},\ldots,\frac{1}{n})$-Bernoulli measures. We note that $\lambda_{n}$ is the only $\sigma_{n}$-invariant measure $\mu$ on $X_{n}$ which pushes forward to $\lambda_{n}^{+}$.

There is also a factor map
$$\theta_{n} \colon X_{n}^{+} \to \T$$
$$\theta_{n}((x_{k})) = \sum_{k=0}^{\infty}\frac{x_{k}}{n^{k+1}}$$
satisfying $\theta_{n}\sigma_{n}^{+} = T_{n}\theta_{n}$. Both $P_{n}$ and $\theta_{n}$ are Lipschitz.

For $p$ dividing $n$ we define the symbolic version of $T_{p}$ by
$$\tau_{p} \colon X_{n}^{+} \to X_{n}^{+}$$
$$\tau_{p}(x)_{m} = \left( p \cdot x_{m} + \left\lfloor \frac{p \cdot x_{m+1}}{n} \right\rfloor \right) \textrm{ mod } n.$$

Then $T_{p} \theta_{n} = \theta_{n}\tau_{p}$ as well, and $\theta_{n}$ is a factor from $\tau_{p}$ onto $T_{p}$. We note that for $pq = n$ we have $\tau_{p}\tau_{q} = \tau_{q}\tau_{p} = \sigma_{n}^{+}$; we occasionally use the fact that if a measure is jointly $\tau_{p},\tau_{q}$-invariant, then it is $\sigma_{n}^{+}$-invariant as well. Both $\tau_{p}, \tau_{q}$ also define homeomorphisms on $X_{n}$ satisfying $P_{n}\tau_{k} = \tau_{k}P_{n}$ for each $k=p,q$, and $\tau_{p}\tau_{q}=\tau_{q}\tau_{p}=\sigma_{n}$ on $X_{n}$.

We will make use of the following easy result: if $Y \subset X_{n}^{+}$ is compact then $\dimH \theta_{n} Y = \frac{\dimH Y}{\log n}$ (the $\log n$ appears in the denominator due to our choice of using base $2$ in our metric on $X_{n}^{+}$). It is also straightforward to check that $\lambda_{n}^{+}$ is the only measure on $X_{n}^{+}$ which pushes forward to $\lambda_{\T}$.

We will also use the following well-known fact: If $\varphi \colon X \to Y$ is a factor map for some $\mathbb{N}^{2}$-actions on compact metric spaces $X$ and $Y$ and $\mu$ is a Borel probability measure on $Y$ invariant under the $\mathbb{N}^{2}$-action, then there exists a Borel probability measure $\mu^{\prime}$, invariant for the action on $X$, such that $\varphi_{*}(\mu^{\prime}) = \mu$.

\subsection{Space of compact subsets}
For a space $X$ we let $\CK(X)$ denote the space of nonempty compact subsets of $X$ equipped with the Hausdorff metric. Given a continuous map $f \colon X \to Y$ we typically write $\tilde{f} \colon \CK(X) \to \CK(Y)$ for the induced map defined by $\tilde{f}(A) = f(A)$. If $f$ is surjective then $\tilde{f}$ is as well, and if $f$ is $K$-Lipschitz then $\tilde{f}$ is as well. 

\subsection{Dimension}
Throughout, $\dimH$ denotes the Hausdorff dimension of a set. We recall some facts about dimension. Let $\mu$ be a Borel measure on some space $X$. The lower Hausdorff dimension of $\mu$ is 
$$\dim_{*}(\mu) = \inf \{\dimH(E) : \mu(E) > 0\}.$$

The lower local dimension of $\mu$ at $x$ is
$$\underline{\dim}(\mu,x) = \liminf_{r \to 0}\frac{\log \mu(B_{r}(x))}{\log r}$$
and a standard fact (see~\cite[Ch. 10]{FalconerBook}) is that
$$\dim_{*}(\mu) = \underset{x \sim \mu}{\textrm{ess} \hspace{.02in}\textrm{inf}} \, \underline{\dim}(\mu,x).$$

The following theorem plays a key role for us.
\begin{theorem}[Hochman-Shmerkin~{\cite[Theorem 10.4]{HochmanShmerkin}}]\label{thm:HochmanShmerkin}
Let $p,q \ge 2$ be multiplicatively independent natural numbers and suppose $\mu$ is a Borel probability measure on $\T$ which is invariant under both $T_{p}$ and $T_{q}$. If $\dim_{*}(\mu)>0$ then $\mu=\lambda_{\T}$.
\end{theorem}

\section{Proof of Theorem~\ref{thm:mainthm}}
We now prove Theorem~\ref{thm:mainthm}. As mentioned earlier, we work symbolically, and the goal of the section is in fact to prove the following symbolic version of Theorem~\ref{thm:mainthm}.
\begin{theorem}\label{thm:symbolicmainthm}
Let $p,q \ge 2$ be multiplicatively independent natural numbers and let $n=pq$. Suppose $\mu \in \CM(\CK(X_{n}^{+}))$ is invariant under both $\tilde{\tau}_{p},\tilde{\tau}_{q}$ and $\dimH Y > 0$ for $\mu$-a.e. $Y$. Then $\mu=\delta_{X_{n}^{+}}$.
\end{theorem}

Before proving Theorem~\ref{thm:symbolicmainthm}, we first show how to deduce Theorem~\ref{thm:mainthm} from Theorem~\ref{thm:symbolicmainthm}. The details are routine; we include them for completeness.

\begin{proposition}\label{prop:liftingIRCs}
Let $2 \le p,q \in \mathbb{N}$ and $n=pq$.
Suppose $\mu$ is an IRC for $T_{p},T_{q}$ acting on $\T$ and $\dimH(Y) > 0$ for $\mu$-a.e. $Y$. Then there exists an IRC $\mu^{\prime}$ for $\tau_{p},\tau_{q}$ acting on $X_{n}^{+}$ such that $(\tilde{\theta}_{n})_{*}(\mu^{\prime}) = \mu$ and $\dimH(Z) > 0$ for $\mu^{\prime}$-a.e. $Z$.
\end{proposition}
\begin{proof}
Let $\mu \in \CM(\CK(\T))$ be as such. The map $\tilde{\theta}_{n}$ semiconjugates the actions of $\tilde{\tau}_{p},\tilde{\tau}_{q}$ with that of $\tilde{T}_{p},\tilde{T}_{q}$, so there exists $\mu^{\prime} \in \CM(\CK(X_{n}^{+}))$ which is $\tilde{\tau}_{p},\tilde{\tau}_{q}$-invariant and satisfies $(\tilde{\theta}_{n})_{*}\mu^{\prime}=\mu$. Since $\theta_{n}$ is Lipschitz, given $Z \in \CK(X_{n}^{+})$ we have
$$\dimH(\theta_{n}(Z)) \le \dimH(Z).$$
Thus
$$\tilde{\theta}_{n}^{-1}\left(\{Y \in \CK(\T) : \dimH(Y) > 0\}\right) \subset \{Z \in \CK(X_{n}^{+}) : \dimH(Z) > 0\}$$
and hence
$$\mu^{\prime}\left(\{Z \in \CK(X_{n}^{+}) : \dimH(Z) > 0\}\right) \ge (\tilde{\theta}_{n})_{*}\mu^{\prime}\left(\{Y \in \CK(T) : \dimH(Y) > 0\}\right)$$
$$= \mu\left(\{Y \in \CK(T) : \dimH(Y) > 0\}\right)= 1.$$
\end{proof}

Now to see how Theorem~\ref{thm:mainthm} follows from Theorem~\ref{thm:symbolicmainthm}, suppose Theorem~\ref{thm:symbolicmainthm} holds, and suppose $\mu \in \CM(\CK(\T))$ is invariant under both $\tilde{T}_{p},\tilde{T}_{q}$ and $\dimH Y > 0$ for $\mu$-a.e. $Y$. By Proposition~\ref{prop:liftingIRCs}, there exists $\mu^{\prime} \in \CM(\CK(X_{n}^{+}))$ which is $\tilde{\tau}_{p},\tilde{\tau}_{q}$-invariant, satisfies $\dimH Y > 0$ for $\mu^{\prime}$-a.e. $Y$, and $(\tilde{\theta}_{n})_{*}\mu^{\prime} = \mu$. Then Theorem~\ref{thm:symbolicmainthm} implies $\mu^{\prime} = \delta_{X_{n}^{+}}$ and hence $\mu = (\tilde{\theta}_{n})_{*}\delta_{X_{n}^{+}} = \delta_{\T}$.

\subsection{Proof of Theorem~\ref{thm:symbolicmainthm}}
We now prove Theorem~\ref{thm:symbolicmainthm}. Fix multiplicatively independent natural numbers $p,q \ge 2$ and let $n=pq$. Suppose $\mu \in \CM(\CK(X_{n}^{+}))$ is invariant under the induced maps of $\tau_{p},\tau_{q}$ and that $\dimH Y > 0$ for $\mu$-a.e. $Y$. We first note that we may assume $\mu$ is jointly ergodic. Indeed, the subset $\{Y \in \CK(X_{n}^{+}) : \dimH Y > 0\}$ is Borel, so almost-every ergodic component in the ergodic decomposition of $\mu$ must give full measure to the set $\{Y \in \CK(X_{n}^{+}) : \dimH Y > 0\}$. If the theorem holds in the ergodic case, then almost-every ergodic component of $\mu$ must be $\delta_{X_{n}^{+}}$, and hence $\mu = \delta_{X_{n}^{+}}$ as well. Thus from here on we assume $\mu$ is jointly ergodic. We define
$$\mathcal{X}^{+} = \{(Y,x) \in K(X_{n}^{+}) \times X_{n}^{+}: x \in Y\}$$
considered as a subspace of $\CK(X_{n}^{+}) \times X_{n}^{+}$. We denote projection onto the first coordinate by
$$\pi_{\CX^{+}} \colon \CX^{+} \to \CK(X_{n}^{+})$$
$$\pi_{\CX^{+}}(Y,x) = Y.$$

For $p$ dividing $n$ we define
$$S_{p} \colon \CX^{+} \to \CX^{+}$$
$$S_{p}(Y,x) = (\tilde{\tau}_{p}(Y),\tau_{p}(x)).$$

Note that $\pi_{\mathcal{X}^{+}}$ is a factor map intertwining $S_{p}$ and $\tilde{\tau}_{p}$, and the fiber over any $Y$ can be identified with $Y$.

We recall (see~\cite{Bowenentropy1}) for a continuous $T \colon X \to X$ and compact $K \subset X$ the topological entropy of the set $K$ is
\begin{equation}\label{eqn:relentropydef}
h_{\textrm{top}}(T,K) = \lim_{\varepsilon \to 0}\limsup_{n \to \infty} \frac{1}{n} \log s_{n}(T,K,\varepsilon)
\end{equation}
where $s_{n}(T,K,\varepsilon)$ is the maximum cardinality of an $\varepsilon$-separated subset of $K$ for the metric $d^{(n)}(x,y) = \max\{d(T^{i}(x),T^{i}(y) : 0 \le i < n\}$

\begin{lemma}\label{lemma:relentropyboundbydim}
For every $Y \in \CK(X_{n}^{+})$, we have $\htop(S_{p}, \pi_{\CX^{+}}^{-1}(Y)) \ge \frac{\log(p)}{\log n} \dimH Y$.
\end{lemma}
\begin{proof}
First note that, letting $d_{\CX^{+}}$ denote the max product metric on $\CX^{+}$, we have that $d_{\CX^{+}}(S_{p}^{j}(Y,x),S_{p}^{j}(Y,y)) = d(\tau_{p}^{j}(x),\tau_{p}^{j}(y))$ for every $j \ge 0$, and hence $\htop(S_{p},\pi_{\mathcal{X}^{+}}^{-1}(Y)) = \htop(\tau_{p},Y)$. Setting $\theta_{n}(Y)=Z \subset \T$, since $\htop(\tau_{p},Y) \ge \htop(T_{p},Z)$ and $\dimH Z = \frac{\dimH Y}{\log n}$, it suffices to prove that $\htop(T_{p},Z) \ge \log(p)\dimH Z$.

Let $s < \dimH Z$. The upper box dimension is bounded below by Hausdorff dimension, so for infinitely many $k$ there exists a $p^{-k}$-separated set $E_{k} \subset Z$ satisfying $|E_{k}| \ge p^{sk}$. For distinct $x,y \in E_{k}$, there exists $0 \le j < k$ such that $1/p^{2} \le p^{j}d_{\T}(x,y) \le 1/p$, and hence $E_{k}$ is a $\frac{1}{p^{2}}$-separated set in $Z$ for the $k$-Bowen metric defined by $T_{p}$. Since the $\limsup$ term in~\eqref{eqn:relentropydef} is nondecreasing as $\varepsilon \to 0$, we have $\htop(T_{p},Z) \ge \log(p) s$. Since $s < \dimH Z$ was arbitrary, this proves the lemma.
\end{proof}

We set $\mathcal{A}^{+} = \pi_{\CX^{+}}^{-1}(\mathcal{B}(\CK(X_{n}^{+})))$
where $\mathcal{B}(\CK(X_{n}^{+}))$ is the Borel $\sigma$-algebra of $\CK(X_{n}^{+})$. The Ledrappier-Walters relative variational principle~\cite[Thm 2.1]{LedrappierWalters} implies
$$\sup\{h_{\rho}(S_{p} \mid \mathcal{A}^{+}) : (S_{p})_{*}\rho = \rho, (\pi_{\CX^{+}})_{*}\rho=\mu\} = \int_{\CK(X_{n}^{+})}\htop(S_{p},\pi_{\CX^{+}}^{-1}(Y))\, d\mu(Y).$$

Since the right hand side of the above is positive by Lemma~\ref{lemma:relentropyboundbydim} and the dimension assumption on $\mu$, there exists some $S_{p}$-invariant measure $\rho$ such that
\begin{equation}\label{eqn:posrelentropy1}
(\pi_{\CX^{+}})_{*}(\rho)=\mu \, \textrm{ and } h_{\rho}(S_{p} \mid \mathcal{A}^{+}) > 0.
\end{equation}

We want to pass to a setting where $\tau_{p},\tau_{q}$ and $S_{p},S_{q}$ act invertibly. One way would be via the joint natural extension; however, with a small amount of work we can use the two-sided shift $X_{n}$, on which $\tau_{p}$ and $\tau_{q}$ already act invertibly, and we opt for this. Similar to before, we define
$$\CX = \{(Y,x) \in \CK(X_{n}) \times X_{n} : x \in Y\} \subset \CK(X_{n}) \times X_{n}.$$

We similarly denote projection onto the first coordinate by
$$\pi_{\CX} \colon \CX \to \CK(X_{n})$$
$$\pi_{\CX}(Y,x) = Y$$
and for $k=p,q$ define
$$G_{k} \colon \CX \to \CX$$
$$G_{k}(Y,x) = \left(\tilde{\tau}_{k}(Y),\tau_{k}(x)\right).$$
Note that each $G_{k}$ is a homeomorphism.
We define
$$\mathcal{A} = \pi_{\CX}^{-1}(\mathcal{B}(\CK(X_{n})))$$
where $\mathcal{B}(\CK(X_{n}))$ is the Borel $\sigma$-algebra of $\CK(X_{n})$. 
There is also a projection map
$$Q \colon \CX \to \CX^{+}$$
$$Q(Y,x) = (\tilde{P}_{n}(Y),P_{n}(x))$$
which is a factor map between $G_{k}$ and $S_{k}$ for $k=p,q$.

Our goal is to lift $\rho$ to a measure $\eta$ on $\CX$ which is $G_{p}$ and $G_{q}$-invariant and whose entropy for $G_{p}$ conditioned on $\mathcal{A}$ is still positive.

\begin{lemma}\label{lemma:gettingeta}
There exists $\eta \in \mathcal{M}(\CX)$ which is jointly $G_{p},G_{q}$\-invariant, satisfies $(\tilde{P}\pi_{\CX})_{*}\eta = \mu$, and $h_{\eta}(G_{p} \mid \mathcal{A}) \coloneq e_{0} > 0$.
\end{lemma}
\begin{proof}
Let $\hat{\CX}^{+}$ denote the natural extension of $S_{p} \colon \CX^{+} \to \CX^{+}$ and define
$$\Phi \colon \CX \to \hat{\CX}^{+}$$
$$\Phi(Y,x) = \left( Q(Y,x), Q(G_{p}^{-1}(Y,x)), Q(G_{p}^{-2}(Y,x)),\ldots \right).$$
This is well-defined since $Q$ is a factor map between $G_{p}$ and $S_{p}$. Then $\Phi G_{p} = \hat{S}_{p} \Phi$, and letting $p_{0} \colon \hat{\CX}^{+} \to \CX^{+}$ denote projection onto the first coordinate, we also have $p_{0} \Phi = Q$.

We claim $\Phi$ is onto. Given a point $\omega = \left((Y^{(0)},x^{(0)}),(Y^{(-1)},x^{(-1)}),\ldots \right) \in \hat{\CX}^{+}$ define
$$Y(\omega) = \{y \in X_{n} : P_{n}\tau_{p}^{-m}(y) \in Y^{(-m)} \textrm{ for all } m \ge 0\}.$$
Choose $z \in X_{n}$ such that $P_{n}\tau_{p}^{-i}(z) = x^{(-i)}$ for every $i$; such $z$ exists since the collection $U_{I} = \{z \in X_{n} : P_{n}\tau_{p}^{-i}(z) = x^{(-i)} \textrm{ for all } 0 \le i \le I\}$ is nested and each is nonempty and compact. Note that such a $z$ belongs to $Y(\omega)$.

We claim
\begin{equation}\label{eqn:Ymequation}
P_{n}(\tau_{p}^{-m}(Y(\omega))) = Y^{(-m)}
\end{equation}
for every $m \ge 0$. Indeed, we clearly have $P_{n}\tau_{p}^{-m}(Y(\omega)) \subset Y^{(-m)}$. If $y \in Y^{(-m)}$, then since $\tau_{p}(Y^{(-i-1)}) = Y^{(-i)}$ for every $i$, we can choose a sequence $y_{i} \in Y^{(-i)}$ with $\tau_{p}(y_{i+1})=y_{i}$ and $y_{m} = y$. As before we can find $w \in X_{n}$ such that $P_{n}\tau_{p}^{-i}(w) = y_{i}$ for every $i$. Then $w \in Y(\omega)$ and $P_{n}\tau_{p}^{-m}(w)=y_{m}=y$, so $Y^{(-m)} \subset P_{n}\tau_{p}^{-m}(Y(\omega))$ as well.
It follows that $\Phi(Y(\omega),z) = \omega$, and hence $\Phi$ is onto.

Let $\hat{\rho}$ denote the natural extension of the measure $\rho $ on $\hat{\CX}^{+}$.

A quick check shows that
\begin{equation}\label{eqn:saturations}
\tau_{p}(Y(\omega)) = Y(\hat{S}_{p}(\omega)).
\end{equation}
We define
$$C(\omega) = \{y \in X_{n} : P_{n}(\tau_{p}^{-m}(y)) = x^{(-m)} \textrm{ for every } m \ge 0\}.$$
Since $x^{(-m)} \in Y^{(-m)}$ for every $m$ we have $C(\omega) \subset Y(\omega)$. Let
$$\CY = \{(\omega,y) \in \hat{\CX}^{+} \times X_{n} : y \in C(\omega)\}.$$
Then $\CY$ is compact and $H \colon \CY \to \CY$ defined by $H(\omega,y) = (\hat{S}_{p}(\omega),\tau_{p}(y))$ is a homeomorphism. Write $\pi_{\CY} \colon \CY \to \hat{\CX}^{+}$ for projection onto the first coordinate, which is a factor map between $H$ and $\hat{S}_{p}$. Then we can find a probability measure $\mu_{0}$ on $\CY$ which is $H$ invariant and whose pushforward under projection onto the first coordinate is $\hat{\rho}$. Define
$$\Psi \colon \CY \to \CX$$
$$\Psi(\omega,y) = (Y(\omega),y).$$
Then $\Psi$ is Borel, and using~\eqref{eqn:saturations} we have that $\Psi H = G_{p} \Psi$.

Define
$$\eta_{0} = \Psi_{*}(\mu_{0}) \in \CM(\CX).$$
Then $\eta_{0}$ is $G_{p}$-invariant. Recall $\mathcal{A} = \pi_{\mathcal{X}}^{-1}(\mathcal{B}(\CK(X_{n})))$. Let $\hat{\mathcal{A}}$ denote the $\sigma$-algebra generated by the first coordinate of $\hat{\CX}^{+}$ (in other words, the preimage of the Borel $\sigma$-algebra on the natural extension of $\tilde{\tau}_{p}$ on $\CK(X_{n}^{+})$). 

We claim $\pi_{\mathcal{Y}}^{-1}(\hat{\mathcal{A}}) = \Psi^{-1}(\mathcal{A})$. Indeed, given how $\left(Y^{(0)},Y^{(-1)},\ldots\right)$ determines $Y(\omega)$, it is easy to see that $\Psi^{-1}(\mathcal{A}) \subset \pi_{\mathcal{Y}}^{-1}(\hat{\mathcal{A}})$. Also, by~\eqref{eqn:Ymequation}, for each $m \ge 0$ the map $\omega = \left((Y^{(0)},x^{(0)}),(Y^{(-1)},x^{(-1)}),\ldots \right) \mapsto Y^{(-m)}$ is a continuous function of $Y(\omega)$, from which it follows that $\pi_{\mathcal{Y}}^{-1}(\hat{\mathcal{A}}) \subset \Psi^{-1}(\mathcal{A})$.

Let $\CE$ be the partition of $\CX$ determined by the zero coordinate in the $X_{n}$-component of $\CX$, and let $\CD$ be the partition of $\hat{\CX}^{+}$ determined by the zero coordinate of $x^{(0)}$. Since $(\omega,y) \in \mathcal{Y}$ implies $y \in C(\omega)$ and hence $P(y) = x^{(0)}$, we have $\Psi^{-1}(\CE) = \pi_{\CY}^{-1}(\CD)$. Moreover, $\Psi$ and $\pi_{\CY}$ intertwine the respective dynamics, so
$$\Psi^{-1}\left(\bigvee_{i=0}^{N-1}G_{p}^{-i}\CE\right) = \pi_{\CY}^{-1}\left(\bigvee_{i=0}^{N-1}\hat{S}_{p}^{-i}\CD\right).$$

Then
$$H_{\eta_{0}}\left(\bigvee_{i=0}^{N-1}G_{p}^{-i}\CE \mid \mathcal{A}\right) = H_{\mu_{0}}\left(\pi_{\CY}^{-1} \bigvee_{i=0}^{N-1}\hat{S}_{p}^{-i}\CD \mid \pi_{\mathcal{Y}}^{-1}\hat{\mathcal{A}}\right)$$
$$=H_{\hat{\rho}}\left(\bigvee_{i=0}^{N-1}\hat{S}_{p}^{-i}\CD \mid \hat{\mathcal{A}}\right)$$
so dividing by $N$ and taking $N \to \infty$ yields
$$h_{\eta_{0}}(G_{p},\CE \mid \mathcal{A}) = h_{\hat{\rho}}(\hat{S}_{p},\CD \mid \hat{\mathcal{A}}).$$

Furthermore, it is not hard to check that both $\CE$ and $\CD$ are conditional entropy generators for $G_{p}, \mathcal{A}$ and $\hat{S}_{p}, \hat{\mathcal{A}}$ respectively, so $h_{\eta_{0}}(G_{p} \mid \mathcal{A}) = h_{\eta_{0}}(G_{p},\CE \mid \mathcal{A})$ and $h_{\hat{\rho}}(\hat{S}_{p} \mid \hat{\mathcal{A}}) = h_{\hat{\rho}}(\hat{S}_{p},\CD \mid \hat{\mathcal{A}})$.

As conditional entropy is preserved when passing to the natural extension, we get that
$$h_{\rho}(S_{p} \mid \mathcal{A}^{+}) = h_{\hat{\rho}}(\hat{S}_{p} \mid \hat{\mathcal{A}}) = h_{\eta_{0}}(G_{p} \mid \mathcal{A}).$$

Note that if $\omega = \left((Y^{(0)},x^{(0)}),(Y^{(-1)},x^{(-1)}),\ldots \right)$ then $Q\Psi(\omega,y) = (Y^{(0)},x^{(0)})$, since $P(Y(\omega)) = Y^{(0)}$ follows from~\eqref{eqn:Ymequation}, and $P(y) = x^{(0)}$ since $y \in C(\omega)$. Thus $Q\Psi = p_{0}\pi_{\mathcal{Y}}$.

Then by~\eqref{eqn:posrelentropy1}, we have $\eta_{0} \in \CM(\CX)$ which is $G_{p}$-invariant and satisfies
$$h_{\eta_{0}}(G_{p} \mid \mathcal{A}) \coloneq e_{0} > 0$$
$$Q_{*}\eta_{0} = Q_{*}\Psi_{*}\mu_{0} = (p_{0}\pi_{\mathcal{Y}})_{*}\mu_{0} = \rho.$$

We now promote $\eta_{0}$ to a measure which is also $G_{q}$-invariant in the usual way. Since $\tau_{q}$ acts invertibly on $X_{n}$, the map $G_{q}$ is invertible on $\CX$. Moreover, $G_{q}$ preserves the $\sigma$-algebra $\mathcal{A}$, so
\begin{equation*}
h_{(G_{q})_{*}\eta_{0}}(G_{p} \mid \mathcal{A}) = h_{\eta_{0}}(G_{p} \mid \mathcal{A}).
\end{equation*}
Define
$$\eta_{m} = \frac{1}{m}\sum_{i=0}^{m-1}(G_{q}^{i})_{*}\eta_{0}.$$
For every $m \ge 1$ we have
$$Q_{*}\eta_{m} = \frac{1}{m}\sum_{i=0}^{m-1}(S_{q}^{i})_{*}\rho$$
and since $\mu$ is $\tilde{\tau}_{q}$-invariant and $(\pi_{\CX^{+}})_{*}\rho = \mu$, this implies
$$(\tilde{P}\pi_{\CX})_{*}\eta_{m} = (\pi_{\CX^{+}}Q)_{*}\eta_{m} = \frac{1}{m}\sum_{i=0}^{m-1}(\tilde{\tau}_{q}^{i})_{*}\mu = \mu.$$

Since conditional entropy is affine (see e.g.~\cite[Sec. 2.6]{DownarowiczBook}), we also have $h_{\eta_{m}}(G_{p} \mid \mathcal{A}) = e_{0}$. Let $\eta$ be a weak* accumulation point of $\eta_{m}$. Then $\eta$ is jointly $G_{p},G_{q}$-invariant, satisfies $(\tilde{P}\pi_{\CX})_{*}\eta = \mu$, and $h_{\eta}(G_{p} \mid \mathcal{A}) \ge e_{0}$, where the last part follows from upper semicontinuity of relative entropy (see either~\cite{MaChen2013} or~\cite[Lemma 3.4]{HuangYeZhang2006}).
\end{proof}

We now have $\eta \in \mathcal{M}(\mathcal{X})$ satisfying Lemma~\ref{lemma:gettingeta}. We can safely assume $\eta$ is ergodic as follows. Decompose $\eta$ into ergodic components for the $\mathbb{Z}^{2}$-action of $G_{p},G_{q}$ on $\CX$. Since $\mu$ is ergodic, almost every ergodic component projects to $\mu$, and since conditional entropy is strongly affine on measures, we may choose such a component $\eta_{\textrm{e}}$ such that $h_{\eta_{\textrm{e}}}(G_{p} \mid \mathcal{A}) > 0$. Thus from here on, we will assume $\eta$ is jointly $G_{p},G_{q}$-ergodic and continue to write $0 < e_{0} = h_{\eta_{\textrm{e}}}(G_{p} \mid \mathcal{A})$.

We define
$$\tilde{\mu} = (\pi_{\mathcal{X}})_{*}\eta \in \CM(\CK(X_{n}))$$
and note that
$$\tilde{P}_{*}\tilde{\mu} = \mu$$
and $\tilde{\mu}$ is jointly ergodic for $\tilde{\tau}_{p},\tilde{\tau}_{q}$. 
It suffices to show that $\tilde{\mu} = \delta_{X_{n}}$.

The next step is to disintegrate $\eta$ over $\tilde{\mu}$, writing
$$\eta = \int_{\CK(X_{n})}\overline{\eta}_{Y}\, d\tilde{\mu}.$$
We define
$$\eta_{Y} = (\pi_{2})_{*}\overline{\eta}_{Y}$$
where $\pi_{2} \colon \CX \to X_{n}$ denotes projection onto the second coordinate. Note each $\eta_{Y}$ is supported on $Y$. For $k=p,q$ we claim
\begin{equation}\label{eqn:disintregatedeqn}
\eta_{\tilde{\tau}_{k}(Y)} = (\tau_{k})_{*}\eta_{Y} \textrm{ for } \tilde{\mu}\textrm{-a.e. } Y \in \CK(X_{n}).
\end{equation}

Indeed, since $\eta$ is $G_{k} = \tilde{\tau}_{k} \times \tau_{k}$-invariant and $\tilde{\mu}$ is $\tilde{\tau}_{k}$-invariant, it is straightforward to check that the family $\eta^{\prime}_{Y} \coloneq (G_{k})_{*}\overline{\eta}_{\tilde{\tau}_{k}^{-1}(Y)}$ is a disintegration of $\eta$. By uniqueness of disintegration we must have $\overline{\eta}_{Y} = (G_{k})_{*}\overline{\eta}_{\tilde{\tau}_{k}^{-1}(Y)}$ for $\tilde{\mu}$-a.e. $Y$, and hence $\overline{\eta}_{\tilde{\tau}_{k}(Y)} = (G_{k})_{*}\overline{\eta}_{Y}$ for $\tilde{\mu}$-a.e. $Y$, and it follows that~\eqref{eqn:disintregatedeqn} holds.

Recall $h_{\eta}(G_{p} \mid \mathcal{A}) = e_{0}$.

\begin{lemma}\label{lemma:posdimlemma}
$\dim_{*}\eta_{Y} \ge e_{0}$ for $\tilde{\mu}$-a.e. $Y$.
\end{lemma}
\begin{proof}
Recall $\CE$ is the partition of $\CX$ based on the zero-coordinate in the $X_{n}$-component of $\CX$, we have $\mathcal{A} = \pi_{\mathcal{X}}^{-1}(\mathcal{B}(\CK(X_{n})))$, and that $\CE$ is an entropy-generator for the conditional entropy $h_{\xi}(G_{p} \mid \CA)$ for every $G_{p}$-invariant probability measure $\xi$ on $\CX$. Set
$$\CE_{N} = \bigvee_{i=0}^{N-1}G_{p}^{-i}\CE$$
and
$$\mathcal{F}_{N} = \bigvee_{i=0}^{N-1}\tau_{p}^{-i}\mathcal{F}$$
where $\mathcal{F}$ is the partition of $X_{n}$ determined by the 0-coordinate. Then
$$\pi_{2}^{-1}\mathcal{F}_{N}(x) = \mathcal{E}_{N}(Y,x)$$
for every $(Y,x) \in \CX$, and since $\eta_{Y} = (\pi_{2})_{*}\overline{\eta}_{Y}$, it follows that
$$\eta_{Y}\left(\mathcal{F}_{N}(x)\right) = \overline{\eta}_{Y}(\mathcal{E}_{N}(x)).$$

Then by the conditional Shannon-McMillan-Breiman Theorem (see e.g.~\cite[App. B]{DownarowiczBook}) we have
\begin{equation}\label{eqn:SMB}
-\frac{1}{N}\log \eta_{Y}(\mathcal{F}_{N}(x)) = -\frac{1}{N}\log \overline{\eta}_{Y}\left(\mathcal{E}_{N}(Y,x)\right) \to e(Y,x)
\end{equation}
for some measurable function $e$ and $\eta$-a.e. $(Y,x)$; the function $e$ is equal to the conditional entropy $h_{\eta_{e}}(G_{p} \mid \mathcal{A})$ on each ergodic component $\eta_{e}$ of $\eta$ (recall $\eta$ is jointly ergodic, but not necessarily ergodic for $G_{p}$). 

Now $G_{q}$ permutes $G_{p}$ ergodic components, preserves $\mathcal{A}$, and is invertible, so preserves the conditional $G_{p}$-entropy of each, and as $\CE$ is a generator for each ergodic component of $G_{p}$, it follows that $e$ is invariant under $G_{q}$ as well. By joint ergodicity, $e$ must be constant almost everywhere, and its value must be $h_{\eta}(G_{p}, \CE \mid \mathcal{A}) = e_{0}$.

Let $\CW \subset \CX$ be the set of points $(Y,x)$ where we have convergence of the left hand side of~\eqref{eqn:SMB} to $e_{0}$. We change `for $\eta$-a.e. $(Y,x)$' into `for $\tilde{\mu}$-a.e. $Y$ and $\eta_{Y}$-a.e. $x$' as follows. We have $\eta(\CW) = 1$. Writing $\CW_{Y} = \{x : (Y,x) \in \CW\}$ and identifying each $\mathcal{W}_{Y}$ with a subset of the fiber $\pi_{\mathcal{X}}^{-1}(Y)$ of $Y$, disintegrating gives $1 = \int_{\CK(X_{n})}\overline{\eta}_{Y}(\CW_{Y})\, d\tilde{\mu}$ so $\overline{\eta}_{Y}(\CW_{Y}) = 1$ and $\eta_{Y}(\CW_{Y}) = 1$ for $\tilde{\mu}$-a.e. $Y$.

Now let $Y \in \CK(X_{n})$ be such that $\eta_{Y}(\CW_{Y}) = 1$. The goal is to prove that $\dim_{*}\eta_{Y} \ge e_{0}$. Let $E$ be a Borel set with $\eta_{Y}(E)>0$. Let $s <e_{0}$. For $\eta_{Y}$-a.e. $x$ we know that $-\frac{1}{N}\log \eta_{Y}(\CF_{N}(x)) \to e_{0}$. Since $s < e_{0}$, for a set of $x$ of full $\eta_{Y}$-measure, there exists $N_{0}(x)$ such that for $N \ge N_{0}(x)$ we have $-\frac{1}{N}\log \eta_{Y}(\CF_{N}(x)) > s$ and hence $\eta_{Y}(\CF_{N}(x)) < 2^{-sN}$. Note that $N_{0}(x)$ depends on $x$ at the moment. Define
$$E_{m} = \{x : \eta_{Y}(\CF_{N}(x)) \le 2^{-sN} \textrm{ for every } N \ge m\}.$$
Then $E_{1} \subset E_{2} \subset \cdots$ and by the above, we know that
$$\eta_{Y}\left(\bigcup_{i=1}^{\infty}E_{i}\right) = 1.$$
Thus
$$\eta_{Y}\left( E \cap \bigcup_{i=1}^{\infty} E_{i} \right) = \eta_{Y}(E)$$
so
$$\eta_{Y}(E \cap E_{N_{0}}) > 0$$
for some $N_{0}$. Set $E^{\prime} = E \cap E_{N_{0}}$. Then $\eta_{Y}(\CF_{N}(x)) \le 2^{-sN}$ for all $x \in E^{\prime}$ and all $N \ge N_{0}$.

Given $x \in X_{n}$, let $C_{N}(x)$ denote the $2N+1$-cylinder set containing $x$ centered at zero in $X_{n}$, so $C_{N}(x) = \{y \in X_{n} : y_{k} = x_{k} \textrm{ for all } |k| \le N\}$. It is straightforward to check from the definition of $\tau_{p}$ that for $x \in X_{n}$ we have $C_{N}(x) \subset \mathcal{F}_{N}(x)$. It follows there exists a constant $C$ (for instance we can take $C=\frac{1}{\eta_{Y}(E^{\prime}})$) such that for every $N \ge N_{0}$ and every $2N+1$-cylinder set $C_{N}(x)$ centered at zero in $X_{n}$, we have $\eta_{Y}(E^{\prime} \cap C_{N}(x)) \le C2^{-Ns}$. 
Replacing the constant $C$ if necessary to handle the finitely many remaining $N< N_{0}$, the mass distribution principle then implies $\dimH E^{\prime} \ge s$ and hence $\dimH E \ge s$. Since $s < e_{0}$ was arbitrary, the lemma follows.
\end{proof}

Our goal now is to prove that $\eta_{Y} = \lambda_{n}$ for $\tilde{\mu}$-a.e. $Y$. Set
$$\varphi \coloneq \theta_{n}P_{n} \colon X_{n} \to \mathbb{T}$$
and let $\alpha_{Y} = \varphi_{*}\eta_{Y}$. We claim
\begin{equation}\label{eqn:backdowntoT}
\dim_{*}\alpha_{Y} \ge c \coloneq \frac{e_{0}}{\log p} > 0 \quad \textrm{ for } \tilde{\mu}\textrm{-a.e. } Y.
\end{equation}

It's enough to show that the above holds for every $Y$ satisfying~\eqref{eqn:SMB}, so let $Y$ be as such. We note that $\alpha_{Y}$ has no atoms: if $z$ was an atom, then one of the at most two $\theta_{n}$-preimages, say $w$, would be an atom for $(P_{n})_{*}\eta_{Y}$, and then every $\mathcal{F}_{N}$-atom which contains a point from $P_{n}^{-1}(w)$ would have $\eta_{Y}$-measure bounded uniformly from below, contradicting that $Y$ satisfies the limit~\eqref{eqn:SMB}. Set
$$\mathcal{I} = \{\left[\frac{j}{n},\frac{j+1}{n}\right) : 0 \le j \le n-1 \}, \qquad \mathcal{I}_{N} = \bigvee_{i=0}^{N-1}T_{p}^{-i}\mathcal{I}$$
and note that the lengths of the intervals in $\mathcal{I}_{N}$ are $\frac{1}{np^{N-1}}$. Since $\theta_{n}$ is injective off a countable set, $\alpha_{Y}$ has no atoms, and $(\tau_{p}^{i}(x))_{0} = j$ if and only if $T_{p}^{i}(\varphi(x)) \in [\frac{j}{n},\frac{j+1}{n})$, we have $\eta_{Y}(\mathcal{F}_{N}(x)) = \alpha_{Y}(\mathcal{I}_{N}(\varphi(x)))$ for $\eta_{Y}$-a.e. $x$, and hence
\begin{equation}\label{eqn:SMBoncircle}
-\frac{1}{N}\log\alpha_{Y}(\mathcal{I}_{N}(z)) \to e_{0} \quad \textrm{ for } \alpha_{Y}\textrm{-a.e. } z.
\end{equation}

If $0 < t < \frac{e_{0}}{\log p}$, then~\eqref{eqn:SMBoncircle} implies $\alpha_{Y}(\mathcal{I}_{N}(z)) \le 2^{-tN \log p} = p^{-tN}$ for all $N$ sufficiently large. Then an analogous argument as the one given at the end of the proof of Lemma~\ref{lemma:posdimlemma} gives~\eqref{eqn:backdowntoT}.

We now show that
\begin{equation}\label{eqn:alphaYs}
\alpha_{Y} = \lambda_{\mathbb{T}} \textrm{ for } \tilde{\mu}\textrm{-a.e. } Y.
\end{equation}

We follow a well-known method. Define
$$\xi = \int_{\CK(X_{n})} \alpha_{Y} * \check{\alpha}_{Y}\, d\tilde{\mu}$$
where as usual $\check{\alpha}_{Y}(A) = \alpha_{Y}(-A)$. For $k=p,q$, since $\varphi$ intertwines $\tau_{k}$ and $T_{k}$, using~\eqref{eqn:disintregatedeqn} we have
$$\alpha_{\tilde{\tau}_{k}Y} = (T_{k})_{*}\alpha_{Y}.$$

Then using invariance of $\tilde{\mu}$ and the fact that the $T_{k}$ are homomorphisms\footnote{Note that the $\tau_{k}$ are not homomorphisms, which is why we descend to the circle for this step.}, we have
$$(T_{k})_{*}\xi = \int_{\CK(X_{n})} (T_{k})_{*}\alpha_{Y} * (T_{k})_{*}\check{\alpha}_{Y}\, d\tilde{\mu} = \int_{\CK(X_{n})} \alpha_{\tilde{\tau}_{k}Y} * \check{\alpha}_{\tilde{\tau}_{k}Y}\, d\tilde{\mu} = \xi$$
so $\xi$ is jointly $T_{p},T_{q}$-invariant. We claim $\dim_{*}\xi \ge c > 0$. To see this, suppose $\xi(A) > 0$. Then we must have
\begin{equation}\label{eqn:convposmeasure}
\alpha_{Y} * \check{\alpha}_{Y} (A) > 0
\end{equation}
for some set of $Y$s with positive $\tilde{\mu}$-measure. Then we may choose a $Y_{0}$ for which both~\eqref{eqn:convposmeasure} and~\eqref{eqn:backdowntoT} hold. Then
$$0 < \alpha_{Y_{0}} * \check{\alpha}_{Y_{0}} (A) = \int_{\mathbb{T}}\alpha_{Y_{0}}(A+z)\, d\alpha_{Y_{0}}$$
so there must exist some $z \in \mathbb{T}$ such that $\alpha_{Y_{0}}(A+z) > 0$. Then
$$\dimH A = \dimH(A+z) \ge \dim_{*}\alpha_{Y_{0}} \ge c$$
and hence $\dim_{*} \xi \ge c$.

Now we have $\xi \in \CM(\T)$ which is $T_{p},T_{q}$-invariant and $\dim_{*}\xi > 0$. By Theorem~\ref{thm:HochmanShmerkin}, this implies $\xi = \lambda_{\T}$. Letting $\hat{\xi}(m)$ denote the $m$th Fourier coefficient of $\xi$, for $m \ne 0$ we have
$$0 = \hat{\xi}(m) = \int_{\CK(X_{n})}\lvert \hat{\alpha}_{Y}(m) \rvert^{2} \, d\tilde{\mu}$$
so $\hat{\alpha}_{Y}(m) = 0$ for $\tilde{\mu}$-a.e. $Y$. Since there are only countably many $m$ there exists a set $\mathcal{V}$ of full $\tilde{\mu}$-measure such that all Fourier coefficients of $\alpha_{Y}$ for $m \ne 0$ vanish for every $Y \in \mathcal{V}$. By uniqueness of Fourier coefficients, we then have $\alpha_{Y} = \lambda_{\T}$ for $\tilde{\mu}$-a.e. $Y$. Thus
$$(\theta_{n})_{*}(P_{n})_{*}\eta_{Y} = \lambda_{\T}$$
for $\tilde{\mu}$-a.e. $Y$, so
\begin{equation}\label{eqn:etaYstobernoulli}
(P_{n})_{*}\eta_{Y} = \lambda_{n}^{+}
\end{equation}
for $\tilde{\mu}$-a.e. $Y$. Since $\tau_{p}\tau_{q}=\sigma_{n}$, by~\eqref{eqn:disintregatedeqn} there exists a full measure set of $Y$ satisfying $(\sigma_{n}^{-j})_{*}\eta_{Y} = \eta_{\tilde{\sigma_{n}}^{-j}Y}$ for all $j \ge 0$. Thus using~\eqref{eqn:etaYstobernoulli} we have
$$(P_{n})_{*}(\sigma_{n}^{-j})_{*}\eta_{Y} = \lambda_{n}^{+}$$
for every $j \ge 0$ from which it follows that $\eta_{Y} = \lambda_{n}$.

Since $\eta_{Y}$ is supported on $Y$, this implies $\tilde{\mu}$-a.e. $Y$ is equal to $X_{n}$, and hence $\tilde{\mu} = \delta_{X_{n}}$. Finally, since $\tilde{P}_{*}\tilde{\mu} = \mu$, it follows that $\mu = \delta_{X_{n}^{+}}$.

\section{Corollaries}
We now prove Corollary~\ref{cor:corollary1}.
\begin{proof}[Proof of Corollary~\ref{cor:corollary1}]
Suppose $Y$ is compact, nonempty, and $s$-Ahlfors-regular for some $s > 0$. First we claim that if $Z$ is an accumulation point of $\{T_{m}(Y) : m \ge 1\}$ then $\dimH Z \ge s$. By the regularity, there exists a probability measure $\nu$ on $Y$ such that $cr^{s} \le \nu(B_{r}(x)) \le Cr^{s}$ for all $x \in Y$ and $0 < r \le 1$ and some constants $c,C$.  Then there exists a constant $D$ such that for each $m$, there exists a measure $\nu_{m}$ with support contained in $T_{m}(Y)$ such that $\nu_{m}(B(z,r)) \le Dr^{s}$ for every $z$ and $0 < r$; indeed, one can set $E_{m} = Y \cap B_{\frac{1}{10m}}(x)$ where $x$ is any point in $Y$, set $D=\frac{C}{c}20^{s}$, and let $\nu_{m}$ be the pushforward of $\frac{1}{\nu(E_{m})}\nu|_{E_{m}}$ under $T_{m}$. Now suppose $T_{m_{i}}(Y) \to Z$ in $\CK(\T)$. Passing to a subsequence, we may assume $\nu_{m_{i}}$ converges weak* to some $\nu$. Portmanteau implies the support of $\nu$ is contained in $Z$. Finally, for any $z \in \T$ and $r > 0$ we have $\nu(B_{r}(z)) \le \liminf_{i} \nu_{m_{i}}(B_{r}(z)) \le Dr^{s}$, so the mass distribution principle implies $\dimH Z \ge s$.

Now consider the $T_{p},T_{q}$-orbit closure $\mathcal{O} = \overline{\{T_{p}^{a}T_{q}^{b}(Y) : a,b \ge 0\}} \subset \CK(\T)$. By the above, every $Z \in \mathcal{O}$ satisfies $\dimH Z \ge s > 0$. Every accumulation point of the sequence $\eta_{N} = \frac{1}{N^{2}}\sum_{a,b=0}^{N-1}\delta_{T_{p}^{a}T_{q}^{b}Y}$ is jointly $\tilde{T}_{p},\tilde{T}_{q}$-invariant, has support in $\mathcal{O}$, and hence by Theorem~\ref{thm:mainthm}, must be equal to $\delta_{\T}$. This implies $\eta_{N} \to \delta_{\T}$, giving the corollary.
\end{proof}

Another corollary is the following. It can also be deduced from the main result in~\cite{RudolphJohnson}.

\begin{corollary}
Suppose $p,q$ are multiplicatively independent and suppose $Y \subset \T$ is nonempty, compact, $T_{p}$-invariant, and $h_{\textrm{top}}(T_{p}|_{Y}) > 0$. Then for every $\varepsilon > 0$, we have
$$\frac{1}{N} \lvert \{0 \le i \le N-1 : T_{q}^{i}(Y) \textrm{ is } \varepsilon\textrm{-dense in } \T\}\rvert \to 1 \quad \textrm{ as } N \to \infty.$$
\end{corollary}
\begin{proof}
By passing to $\bigcap_{k=0}^{\infty}T_{p}^{k}(Y) \subset Y$ if necessary, we can assume without loss of generality that $T_{p}(Y)=Y$ and $h_{\textrm{top}}(T_{p}|_{Y}) = e > 0$. Since $T_{p}$ and $T_{q}$ commute and $T_{q}$ is finite-to-one, for every $k \ge 0$ we have $T_{p}(T_{q}^{k}(Y)) = T_{q}^{k}(Y)$ and $h_{\textrm{top}}(T_{p}|_{T_{q}^{k}(Y)}) = e$. Consider
$$\mathcal{I}_{p} = \{Z \in \CK(\mathbb{T}) : T_{p}(Z)=Z \textrm{ and } h_{\textrm{top}}(T_{p}|_{Z}) \ge e\}.$$
The invariance property is closed in $\CK(\T)$, and since topological entropy is upper semicontinuous on the space of $T_{p}$-invariant subsets in $\CK(\mathbb{T})$, it follows that $\mathcal{I}_{p}$ is closed, and hence the orbit closure $\overline{\{T_{q}^{k}(Y) : k \ge 1\}}$ is contained in $\mathcal{I}_{p}$. Then every weak* accumulation point of $\mu_{N} = \frac{1}{N}\sum_{k=0}^{N-1}\delta_{T_{q}^{k}(Y)}$ is jointly $\tilde{T}_{p}, \tilde{T}_{q}$-invariant and has support contained in $\mathcal{I}_{p}$. Furthermore, if $Z \in \mathcal{I}_{p}$ then $\dimH Z \ge \frac{e}{\log p}$ (see e.g.~\cite{Furstenberg1967}), so combining this with Theorem~\ref{thm:mainthm}, every weak* accumulation point of $\mu_{N}$ is $\delta_{\T}$, and hence $\mu_{N} \to \delta_{\T}$, which implies the corollary.
\end{proof}

\bibliographystyle{plainurl}
\bibliography{T2T3IRCbib}
\end{document}